\documentclass[12pt]{amsart}

\usepackage{amsthm, amssymb, amsfonts}

\usepackage{mathtext}
\usepackage[cp1251]{inputenc}

\usepackage{url}

\newcommand\R{\mathbb R}

\newcommand{\Cn}{{\mathbb C}^n}

\begin{document}

\title{\bf {Invertibility of quasiconformal operators}}

\author{V.A.~Zorich}

\date{30.04.2021}

\begin{abstract}
The global homeomorphism theorem for quasiconformal maps describes the following specifically higher-dimensional phenomenon:
{\em Locally invertible quasiconformal mapping $f: {\R}^{n} \to {\R}^{n}$ is globally invertible provided
$n > 2$.}

We prove the following operator version of the global homeomorphism theorem.
{\em If the operator $ f: H \to H $ acting in the Hilbert space $ H $ is locally invertible and is an operator of bounded distortion, then it is globally invertible.}
\end{abstract}



\maketitle
\markboth{Vladimir Zorich}{Invertibility of quasiconformal operators}

{\bf 1.~Introduction.}~

\vskip2mm
The global homeomorphism theorem for quasiconformal maps describes the following
specifically higher-dimensional phenomenon (formulated in [1] for $ n = 3 $; proved later in [2]):

{\em Locally invertible quasiconformal mapping
$f: {\R}^{n} \to {\R}^{n}$
is globally invertible provided $n > 2$.}

The theorem has several beautiful generalizations by a number of authors;
they are mentioned in the surveys [3], [4], along with several open problems.
Here we will discuss one of them:
the possibility to extend
the theorem
to nonlinear operators of bounded distortion.

The necessary definitions and facts will be given below.

\vskip5mm
{\bf {2.~Several concepts and facts.}}
\vskip2mm

Let $f:H \to H$ be  a nonlinear smooth operator acting in the Hilbert (or even Banach) space $H$.
Let the tangent operator  $f'(x)$ have a continuous inverse at each point
$x\in H$. This ensures  locally invertibility of the operator $f$.

{\em The coefficient of quasiconformality} of the operator $f$ at the point $x \in H$
is the value
$$
k_{f'(x)}:=\Vert f'(x)\Vert\cdot\Vert (f'(x))^{-1}\Vert.
$$
Geometrically it means
the ratio of the semi-major axis to the semi-minor axis of the ellipsoid that is the image of the unit ball under the linear mapping $f'(x)$.

(A similar value occurs in the operator theory, where it was introduced and used by Banach to measure the distance between different norms in a vector space.{\footnote {Banach considered the product of  norms of the linear operator and its inverse acting between two normed spaces, followed by taking the lower bound of this product over all such linear operators.}})

The quasiconformality coefficient of the mapping $f$ at the point $x$ is usually denoted by the symbol
$k_f(x)$, and its upper bound $ \sup_{x \in D} k_f(x)$ over the domain $D$  of the mapping
$f$ is usually denoted by $ k_f $. The value $ k_f $ is called the {\em coefficient of quasiconformality of the mapping $f$ in the  domain}. In our case $D$ is the whole space $ H $. If $ k_f $ is finite, then the mapping $ f $ is said to have {\em bounded distortion}, or that it is {\em quasiconformal}.
If $ k_f = 1 $, the mapping $ f $ {\em is conformal}.

\vskip2mm
In connection with the theorem on global homeomorphism for quasiconformal mappings of the Euclidean space, mentioned in the introduction, the following question naturally arises: is such a theorem  also true in the infinite-dimensional setting of nonlinear operators of bounded distortion?
An operator analogue of the global homeomorphism theorem would give the following
principle of invertibility:
\vskip2mm
{\em If the operator $ f: H \to H $ acting in the Hilbert space $ H $ is locally invertible and it is an operator of bounded distortion, then it is globally invertible.}
\vskip2mm

Hence, in this case, the equation $ f(x) = y $ has a solution for any right-hand side $ y \in H $, and that solution is unique.
\vskip2mm
It may happens (as it happens in the theory of operators) that properties of the space (Hilbert, Banach) are essential.
We recall that Nevanlinna [5] proved the Liouville theorem on the M\"obius group of conformal mappings even in the case of a Hilbert space. However, I do not know any other results about  just quasiconformal (not quasi-isometric) mappings in the infinite-dimensional case. For quasi-isometric operators (and for non-contractive ones) acting in Banach spaces the following result by John [6] is known: if such an operator is locally invertible, then it is also globally invertible. This statement holds in any dimension (either finite or infinite). A quasi-isometric operator is, of course, a bounded distortion operator; its quasiconformality coefficient is bounded, since the changes of length elements are bounded at any point of the domain, and this estimate is uniform. But not every operator of bounded distortion is a quasi-isometric
operator, as well as not every conformal mapping of a domain on a plane is a quasi-isometric mapping.

\vskip5mm

{\bf {3.~Invertibility of quasiconformal operators.}}
\vskip2mm
Generalizing the concept of conformal (and quasiconformal) mapping of domains of Riemannian manifolds of the same dimension M.Gromov suggested that a mapping of metric spaces is to be considered conformal (respectively, quasiconformal), for example the mapping $ F: \R^{m} \to \R^{n} $ ~ ($ {m} \geq {n} $) if
the image of any infinitely small ball
at each point of the domain
is transformed into an infinitesimal ball (respectively, into an ellipsoid  of uniformly bounded eccentricity) [7].
For instance, any entire holomorphic function $ f: \Cn \to {\mathbb C} $ defines a mapping conformal in the sense of Gromov.

In connection with this extension of the concepts of conformality and quasiconformality,
Gromov naturally wondered what facts of the classical theory may apply to these mappings.
In particular, is it true that {\em if the mapping $ F: \R^{n + 1} \to \R^{n} $
is conformal and bounded,
then for $ n \geq 2 $ it is a constant}?{\footnote {For the information of the reader we note
that on the very first page of the article [7], in the footnote, the author gives the reference where a more complete text of his work in English is presented.

In particular, there the reader can find an extended interpretation of conformity and quasiconformality,
as well as a formulation of the question about the Liouville theorem for such mappings.
This text is now available at
\small {\url {https://www.ihes.fr/~gromov/wp-content/uploads/2018/08/problems-sept2014-copy.pdf}}.
}}

Investigating this problem we came to a construction
which confirmed the validity of such a Liouville-type theorem [8].
Now it seems natural, although it did not immediately become clear that the same construction could work
for investigating the problem of invertibility of operators of bounded distortion.
The point is that mappings that are quasiconformal in the sense of Gromov, in contrast to classical quasiconformal mappings, can act between spaces of different dimensions. This circumstance forced us to expand the classical toolkit of conformal invariants adapting it to such a situation.

As a result it turned out to be possible to prove that {\em if the operator $ f: H \to H $ acting in the Hilbert space $ H $ is locally invertible and an operator of bounded distortion, then the inverse image of any subspace
$ \R^3 \subset H $ of the target space is a non-singular connected three-dimensional surface
$ \tilde S^3 $ properly embedded in the ambient space, while the restriction
$ f|_{\tilde S^3} $ of the original map onto
the surface $ \tilde S^3 $ is quasiconformal, injective, and $ f (\tilde S^3) = \R^3 $.}

\vskip2mm
This statement already implies that {\em the mapping $ f: H \to H $ is globally invertible}.
\vskip2mm
The arguments proving this  auxiliary key statement and its corollary indicated above will be given in the next section.

\vskip5mm
{\bf {4.~Outline of the proof.}}
\vskip2mm

We assume the following normalization by setting $ f (0) = 0 $. Take the germ of the inverse mapping at the point $ 0 \in H $.
Let the germ be defined in some
ball $ B(r_0) $ of radius $ r_0> 0 $.
We will continue this germ along the rays emanating from the point $ 0 $.

In some directions such an extension can go
along the entire ray. Along other directions there may appear singular points that prevent continuation.

Let $ \gamma $ denote the part of the ray from the boundary of the ball $ B (r_0) $, where the original germ of the inverse mapping exits, going to the singular point (that part is finite if continuation along the entire ray is impossible).

We note right away that if the ray along which we continue the original germ of the inverse mapping
abuts against a singular point that prevents further continuation, then
preimage $ \tilde \gamma $ of the specified piece $ \gamma $ of this ray
turns out to be a curve of infinite length, since such a curve cannot lead to any point of the mapped space $ H = \tilde H $. (For further convenience the space-preimage  will be denoted by $ \tilde H $).

We will take in the space of the image $ H $ the subspace
$ \R^2 $ --- a plane, and we will construct
the inverse image of $ \R^2 $ lifting the rays of this subspace
emanating from the point $ 0 $.

The union
$ \Gamma = \{\gamma \} $
of pieces $ \gamma $ of such rays  plus the circle $ \R^2 \cap B(r_0) $ form a star  subdomain $ S^2 $ of the plane $ \R^2 $, and the collection $ \tilde \Gamma = \{\tilde \gamma \} $ of their preimages
plus the preimage of this circle
form a surface $ \tilde S^2 $ without boundary in the space $ \tilde H $.

If we show (this will be done below) that the family $ \tilde \Gamma = \{\tilde \gamma \} $ of curves
$ \tilde \gamma $ on the surface
$ \tilde S^2 $ is a family of curves of conformal module equal to zero ($ M_2 (\tilde \Gamma) = 0 $), then, due to the quasiconformality of the mapping $ f $, the conformal module of the family $ \Gamma = \{\gamma \} $ has to be equal to zero as well ($ M_2 (\Gamma) = 0 $). Then there should be very few singular points preventing such an extension along the rays of the space $ \R^2 $, more precisely, this set should be a set of conformal capacity equal to zero.

Now, following the
proof of the global homeomorphism theorem for locally invertible quasiconformal mappings of the space
$ \R^3 $, one can show that the procedure for such an extension of the germ of the inverse mapping along the rays of any subspace $ \R^3 \subset H $
cannot have singular points at all.\footnote {Having indicated the key idea we do not stop here on the description of the entire technological chain which leads to the conclusion that in the space $ \R^3 $ there will be no singular points at all preventing the procedure of continuation of the germ of the inverse mapping. The details can be found in [2].}

So, it happens that the inverse image of any
subspace $ \R^3 \subset H $ should be
a three-dimensional surface $ \tilde S^3 \subset \tilde H $,
and the restriction of the original mapping $ f $
onto this surface is not only quasiconformal but it is injective and $ f (\tilde S^3) = \R^3 $.

Now, as a corollary, we  conclude
that the initial mapping $ f: H \to H $ is indeed the globally invertible one.

\vskip3mm
It remains to check the equality $ M_2 (\tilde \Gamma) = 0 $.

\vskip2mm

Let us show first that if the mapping $ f: H \to H $
is quasiconformal, then the area of the geodesic circle with center $ 0 $ of radius $ r $ on
the surface $ \tilde S^2 $
increases not faster than $ O(r^2) $
as $ r \to \infty $.

Indeed, since the increment of the area of the circle can be written in the form $ L(r) dr $, where $ L(r) $ is the length of the circle of radius $ r $, then it suffices to check that the length of the circle
increases no faster than $ O(r) $ as $ r \to \infty $.

By assuming the opposite we see that the following situation is inevitable.
Steep waves should appear on this circle, steep in the sense that their height is much greater than their length. And this would contradict the quasiconformality of the mapping $ f $.

This is useful and important observation, so we will reformulate it in a general way
omitting details which are not essential.

\vskip2mm

Take the straight (real) line $ \gamma $ in the Hilbert space $ H $, and let
$ I = [- \varepsilon, \varepsilon] $ be a small segment on $ \gamma $.
Let $ \tilde H $ be another copy of the same Hilbert space $ H $. Take
the curve $ \tilde \gamma $ which differs from the straight line $ \gamma $ only in the way that the segment
$ I $ is replaced by some curve (splash) $ \tilde I $.
Let $ f: \tilde H \to H $ be a quasiconformal mapping such that the curve $ \tilde \gamma $ goes over
into the straight line $ \gamma $ and the surge (wave) $ \tilde I $ goes into the segment $ I $.
The statement is that if the diameter $ \tilde I $ is much larger than the diameter $ I $ (i.e. the height or the amplitude of the wave is much larger than its length), then the coefficient of quasiconformality of the mapping $ f $ must be large.
Moreover, the greater the specified ratio is,
the greater will be the coefficient of quasiconformality.

\vskip2mm

Let us clarify
the assertion by an illustrative example.

Take two copies of the standard Euclidean plane $ \R^{2} $ and $ \tilde \R^{2} $ with Cartesian coordinates $ (x, y) $ and $ (\tilde x, \tilde y) $
respectively. In the $ \R^{2} $ plane mark the straight line $ \gamma $ --- the $ x $ axis, and in the
$ \tilde \R^{2} $ plane instead of the $ \tilde x $ axis consider the curve $ \tilde \gamma $
that coincides with the $ \tilde x $ axis outside a small neighborhood of the origin
and a small segment of this axis containing the origin is replaced by a deep loop (pit, wave) of a completely arbitrary shape.
Suppose we have a quasiconformal mapping $ f: \tilde \R^{2} \to \R^{2} $ such that
the curve $ \tilde \gamma $ transforms to the straight line $ \gamma $. We claim that the coefficient of
quasiconformality of the mapping $ f: \tilde \R^{2} \to \R^{2} $ must be very large, and it will be larger the larger will be depth of the well compared to its entrance (the greater will be the amplitude of the wave  compared to its length).

Note that each of the two regions into which the curve $ \tilde \gamma $ divides the plane $ \tilde \R^{2} $, admits even a conformal mapping onto its upper or lower half-plane, onto which the axis
$ \gamma $ divides the plane $ \R^{2} $.
But we are talking on
a mapping transforming
the domain of the plane $ \tilde \R^{2} $, containing the loop, into a neighborhood of the origin of the plane $ \R^{2} $, and the map straightens
a relatively deep well (or a wave of the amplitude large in comparison with
its length).

This statement can be proved by comparing the conformal moduli of suitable families of curves. For example,
in the $ \R^{2} $ plane take the above-mentioned small neighborhood of the origin on the $ x $ axis
(i.e., on the curve $ \gamma $), take on the same axis $ x $ the exterior of the unit neighborhood of the origin, and
consider
the family of curves connecting these two sets in the $ \R^{2} $ plane.
The conformal modulus of such a family is the same as the conformal capacity of the corresponding condenser.
The conformal modulus of this family of curves is small if the considered neighborhood of the origin is small
(since one of the plates of the condenser is relatively small). However, the conformal modulus of the inverse image of this family of curves cannot be small if
the wave amplitude on the $ \tilde \gamma $ curve is large compared to the wavelength.

The same considerations are applicable in the case when a similar type of curve $ \tilde \gamma $ containing a steep wave
straightens under a quasiconformal mapping $ f: \tilde \R^{n} \to \R^{n} $ of a Euclidean space of any finite dimension. It is also clear that the above reasoning is essentially local.
It will also remain valid if, instead of the mapping $ f: \tilde \R^{n} \to \R^{n} $, we consider the mapping
$ f: \tilde S^{n} \to \R^{n} $ of the Riemannian manifold $ \tilde S^{n} $ obtained by the quasiconformal
deformation of $ \tilde \R^{n} $.

\vskip2mm

Now let us return to the family of curves $ \tilde \Gamma $ on the surface $ \tilde S^2 $ and show that its conformal modulus is equal to zero ($ M_2 (\tilde \Gamma) = 0 $).

Indeed, by construction
curves of the family $ \tilde \Gamma $
have infinite length, going to infinity along the surface $ \tilde S^2 $
(i.e. leaving any compact set of $ \tilde S^2 $).
Thus, the function $ \rho = \alpha ~ (r \ln r)^{- 1} $ (where $ r $ is the geodesic distance from a point of the surface $ \tilde S^2 $ to the starting point $ \tilde 0 = 0 $, and
$ \alpha> 0 $) is
an admissible function for this family of curves. Indeed, $ \int_{\tilde \gamma} \rho = \infty $ for any
curve $ \tilde \gamma \in \tilde \Gamma $. On the other hand, the integral of the square of this function
over that part of the surface that lies outside the neighborhood of the point $ \tilde 0 $,
is finite, and tends to zero as $ \alpha \to 0 $. Hence, $ M_2 (\tilde \Gamma) = 0 $.

(Here we used the fact that the area of the geodesic circle of radius $ r $ on the surface
$ \tilde S^2 $ increases
no faster than $ O(r^2) $ as $ r \to \infty $.)

\vskip3mm
Finally, let us make the last
explanatory remark.
If the procedure for constructing the inverse image of the subspace $ \R^3 \subset H $ is done for each such subspace, then, on the one hand, we exhaust the entire space $ H $ of the image. On the other hand, the corresponding preimages $ \tilde S^3 $ all together
form a subdomain $ \tilde D \subset \tilde H $ of the space-preimage  $ \tilde H $, and the restriction
$ f|_{\tilde D} $ of the mapping $ f $
onto the domain $ \tilde D $ is injective by construction and $ f (\tilde D) = H $. This means that in fact
$ \tilde D = \tilde H $.
Otherwise, one could
take a point
$ p \in \tilde H \setminus \tilde D $ and connect it to the point $ \tilde 0 = 0 $ by a segment.
The part of this segment adjacent to the point $ \tilde 0 = 0 $ certainly lies in the domain $ \tilde D $.
So, in this segment there is such a point, that in any its neighborhood there are points of the domain
$ \tilde D $ and points that do not belong to the region $ \tilde D $. But this is not possible if the mapping $ f: \tilde H \to H $ is locally invertible, $ f (\tilde D) = H $ and $ f|_{\tilde D} $ is injective.

\vskip8mm

\centerline {\bf REFERENCES}
\vskip3mm

[{\bf 1}] M.A.~Lavrent'ev, “On a differential criterion for homeomorphic
mappings of
three-dimensional domains”,
Dokl. Acad. Science USSR, 20 (1938), 241-242.

[Lavrentieff, M.
Sur un crit\`ere diff\'erentiel des transformations homeomorphes des domaines \`a trois dimensions. (French)
C.R. (Dokl.) Acad. Sci. URSS, Ser. 20, 241-242 (1938).]


[{\bf 2}] V.A.~Zorich, “A theorem of M.A.Lavrent'ev on quasiconformal space maps”, Math. USSR-Sb., 3:3 (1967), 389–403.



[{\bf 3}] V.A.~Zorich, “The global homeomorphism theorem for space quasiconformal mappings, its developement and related open problems”, Lecture Notes in Math., 1508 (1992), 132-148.

[{\bf 4}]
V.A.~Zorich, “Quasi-conformal maps and the asymptotic geometry of manifolds”, Russian Math. Surveys, 57:3 (2002), 437–462.

[{\bf 5}] R.~Nevanlinna, “On differentiable mappings”,
Analytic functions. Princeton Math. Series 24 (1960), 3-9.



[{\bf 6}] F.~John,
“On quasi-isometric mappings”,
Comm. Pure and Appl. Math., 22 (1969), 41-66.

[{\bf 7}] M.L.~Gromov, “Colorful categories”, Russian Math. Surveys 70: 4 (2015), 591-655.

[{\bf 8}] V.A.~Zorich, “Conformality in the sense of Gromov and a generalized Liouville theorem”,
http://arxiv.org/abs/2108.00945

\vskip5mm

{\small {\bf V.A.~Zorich}

Lomonosov Moscow State University

\vskip1mm
E-mail: {\bf vzor@mccme.ru}}

\end{document}